\documentclass[12pt, reqno]{amsart}
\usepackage{amsmath}
\usepackage{amsthm}
\usepackage{amssymb}
\usepackage[abbrev]{amsrefs}
\usepackage{mathrsfs}
\usepackage[usenames]{color}
\usepackage{ifthen}
\usepackage{yhmath}
\usepackage{graphicx}
\usepackage[all]{xy}
\usepackage{bm}
\usepackage{bbm}
\usepackage{enumitem}
\usepackage{cases}

\makeatletter

\newcommand{\leqnos}{\tagsleft@true\let\veqno\@@leqno}
\newcommand{\reqnos}{\tagsleft@false\let\veqno\@@eqno}
\reqnos
\makeatother

\AtBeginDocument{\def\MR#1{}}
\newtheorem{thm}{}[section]
\newtheorem{pbm}{}
\newtheorem{theorem}[thm]{Theorem}

\newtheorem{proposition}[thm]{Proposition}

\newtheorem{problem}[pbm]{Problem}
\theoremstyle{definition}
\newtheorem{definition}[thm]{Definition}
\theoremstyle{remark}

\newtheorem{question}[thm]{Question}

\numberwithin{equation}{section}
\allowdisplaybreaks
\newcommand{\ceil}[1]{\left\lceil#1\right\rceil}
\newcommand{\abs}[1]{\left\lvert#1\right\rvert}
\newcommand{\norm}[1]{\left\lVert#1\right\rVert}
\newcommand{\enbrace}[1]{\left\lbrace#1\right\rbrace}
\newcommand{\enpar}[1]{\left(#1\right)}
\newcommand{\Ft}{\ensuremath{\mathcal{F}}}

\newcommand{\Rt}{\ensuremath{\mathcal{R}}}

\newcommand{\xx}{\ensuremath{\bm{x}}}
\newcommand{\xg}{\ensuremath{{\bm{g}}}}
\newcommand{\yy}{\ensuremath{\bm{y}}}
\newcommand{\XX}{\ensuremath{\mathbb{X}}}
\newcommand{\YY}{\ensuremath{\mathbb{Y}}}
\newcommand{\XB}{\ensuremath{\mathcal{X}}}
\newcommand{\Ind}{\ensuremath{\mathbbm{1}}}
\newcommand{\GG}{\ensuremath{\mathcal{G}}}
\newcommand{\FF}{\ensuremath{\mathbb{F}}}
\newcommand{\RR}{\ensuremath{\mathbb{R}}}

\newcommand{\NN}{\ensuremath{\mathbb{N}}}
\newcommand{\ZZ}{\ensuremath{\mathbb{Z}}}
\newcommand{\YB}{\ensuremath{\mathcal{Y}}}
\newcommand{\lsdf}{\ensuremath{\bm{\varphi_l}}}
\newcommand{\usdf}{\ensuremath{\bm{\varphi_u}}}
\newcommand{\prim}{\ensuremath{\bm{\sigma}}}
\newcommand{\Cu}{\ensuremath{\mathcal{Q}}}
\newcommand{\Sym}{\ensuremath{\mathbb{S}}}
\newcommand{\pp}{\ensuremath{{\bm{p}}}}
\newcommand{\Ct}{\ensuremath{\mathcal{C}}}
\newcommand{\BB}{\ensuremath{\mathcal{B}}}
\newcommand{\Zt}{\ensuremath{\mathcal{Z}}}
\DeclareMathOperator{\sgn}{sign}
\newcommand{\BV}{\ensuremath{\mathrm{BV}}}

\newcommand{\unc}{\ensuremath{\bm{k}}}

\newcommand{\ww}{\ensuremath{\bm{w}}}

\newcommand{\DD}{\ensuremath{\mathbb{D}}}

\newcommand{\EE}{\ensuremath{\mathbb{E}}}

\DeclareMathOperator{\spn}{span}
\DeclareMathOperator{\supp}{supp}
\def\cA{{\mathcal A}}

\def\cD{{\mathcal D}}

\def\cH{{\mathcal H}}

\def\cR{{\mathcal R}}

\def\cT{{\mathcal T}}

\title{Twenty-five  years of greedy bases}
\begin{document}
%------------------------------------------------------------------------
\author[Albiac]{Fernando Albiac}
\address{Department of Mathematics, Statistics, and Computer Sciencies--InaMat2 \\
Universidad P\'ublica de Navarra\\
Campus de Arrosad\'{i}a\\
Pamplona\\
31006 Spain}
\email{fernando.albiac@unavarra.es}
%------------------------------------------------------------------------
\author[Ansorena]{Jos\'e L. Ansorena}
\address{Department of Mathematics and Computer Sciences\\
Universidad de La Rioja\\
Logro\~no\\
26004 Spain}
\email{joseluis.ansorena@unirioja.es}
%------------------------------------------------------------------------
\author[Temlyakov]{Vladimir Temlyakov}
\address{Deparment of Mathematics\\
University of South Carolina\\ Columbia, SC 29208 USA;\newline
Steklov Mathematical Institute of Russian Academy of Sciences, Moscow, Russia;\newline Lomonosov Moscow State University; \newline Moscow Center of Fundamental and Applied Mathematics.}
\email{temlyak@math.sc.edu}
%------------------------------------------------------------------------
%------------------------------------------------------------------------
%------------------------------------------------------------------------
\begin{abstract} 
Although the basic idea behind the notion of a greedy basis had been around for some time, the formal development of a theory of greedy bases was initiated in 1999 with the publication of the article \cite{KoTe1999} by Konyagin and Temlyakov.  The theoretical simplicity of the thresholding greedy algorithm became a model for a procedure widely used in numerical applications and the subject of greedy bases evolved very rapidly from the point of view of approximation theory.  The idea of studying greedy bases and related greedy algorithms attracted also the attention of researchers  with a  classical Banach space theory background. From the more abstract point of functional analysis, the theory of greedy bases and its derivates evolved very fast as many fundamental results were discovered and new ramifications branched out.  Hundreds of papers on greedy-like bases and several monographs have been written since  the appearance of the foundational paper \cite{KoTe1999}. After twenty five years, the theory  is very much alive  and it continues to be a very active research topic both for functional analysts and for researchers interested in the applied nature of nonlinear approximation alike. This is why  we believe it is a good moment   to gather a selection of 25 open problems (one per year since 1999!) whose solution would contribute to advance the state of art of this beautiful topic. 
\end{abstract}
%------------------------------------------------------------------------
\thanks{The research of F. Albiac and J.L. Ansorena was partially funded by the Spanish Ministry for Science and Innovation under Grant PID2022-138342NB-I00 for \emph{Functional analysis methods in approximation theory and applications}.
The research of V.N. Temlyakov was supported by the Russian Science Foundation under grant no. 23-71-30001, https://rscf.ru/project/23-71-30001/, and performed at Lomonosov Moscow State University.}
%------------------------------------------------------------------------
\maketitle
%------------------------------------------------------------------------
\section{Introduction}\noindent 
%------------------------------------------------------------------------
Greedy algorithms provide sparse representations (or approximations) of a given image/signal in terms of a given system of elements of the ambient space. In a mathematical setting, an image or signal 
is considered to be a function  of a Banach space. For instance, a two-dimensional image can be viewed as a function of two variables belonging to the  Hilbert space $L_2$ or, more generally, to the Banach space $L_p$, $1\le p\le \infty$.  Usually we assume that the system used for representation has some natural properties and we call it  a \emph{dictionary}.  For an element $f$ from a Banach space $\XX$ and a fixed $m$, we consider  approximants which are linear combinations of $m$ terms from a dictionary $\cD$. We call such an approximant an $m$-\emph{term approximant} of $f$ with respect to $\cD$. 

In sparse approximation, a greedy algorithm is an algorithm that uses a {\it greedy step} in searching for a new element to be added to a given $m$-term approximant. By a {\it greedy step}, we mean one 
which maximizes a certain functional determined by information from the previous steps of the algorithm. We obtain different types of greedy algorithms by varying the above-mentioned functional and also by using different ways of constructing (choosing coefficients of the linear combination) the $m$-term approximant from previously selected $m$ elements of the dictionary. 

A classical problem of mathematical and numerical analysis that goes back to the origins of Taylor's and Fourier's expansions, is to approximately represent  a given function from a space. The first step to solve the representation problem is to choose a representation system. Traditionally, a representation system   has some  natural features such as minimality, or orthogonality, i.e., a simple structure which allows nice computational properties. The most typical representation systems are the trigonometric system
\[
x\mapsto e^{ i  kx}, \quad k\in\ZZ,
\] 
the algebraic system 
\[
x\mapsto x^k, \quad k\in\ZZ,\, k\ge 0,
\] 
the spline system, the wavelet system and their multivariate versions. In general we may speak of a basis

 $\XB=(\xx_k)_{k=1}^\infty$
 (in a sense that we will specify below) in a Banach (or quasi-Banach) space $\XX$.

The second step to solve the representation problem is to choose the form of the approximant to be built from the chosen representation system $\XB$. In a classical way which was used for centuries,  given $m\in\NN$,
an approximant $a_m$ is a polynomial  of order $m$
 with respect to $\XB$, 
 that is, 
\[
a_m=\sum_{k=1}^{m} c_k\, \xx_k
\]
for some scalars $c_k$, $k=1$, \dots, $m$.
In numerical analysis and approximation theory it was understood that in many problems from signal/image processing it is more beneficial to use an $m$-term approximant with respect to $\XB$ than a polynomial of order $m$.   This means that for $f\in \XX$ we look for an approximant of the form 
\[
a_m(f):=\sum_{k\in\Lambda(f)}c_k\, \xx_k,
\]
where $\Lambda(f)$ is a set of $m$ indices which is determined by $f$.

The third step to solve the representation problem is to choose a method of construction of the approximant. In linear theory, partial sums of the corresponding expansion of $f$ with respect to the basis $\XB$ is a standard method. It turns out that greedy approximants are natural substitutes for the partial sums in nonlinear theory. 

We emphasize that, although in the beginning this theory developed  within the framework of approximation theory, soon after the appearance of \cite{KoTe1999} the center of activities in the area moved to functional analysis. Indeed, the introduction of new types of bases and the achievement of their characterizations in terms of classical properties from Banach space theory caught the attention of the specialists, who gave impetus to the theory and set the foundations for a fruitful and novel research topic. 
As a result we now have different notations 
for the same objects, which come from approximation theory (see e.g. \cite{VTbook} and \cite{VTsparse}) and functional analysis.  In keeping with current usage, in the next section we will present the notation and terminology in the way it is nowadays used in the modern functional analysis approach to the subject, as the reader  can find in \cite{AlbiacKalton2016}*{Chapter 10} or \cite{AABW2021}. After the preliminary Section~\ref{NotTerm}, we present a selection of topics that reflect the state of art of the theory and suggest within each section the problems  that we believe should be addressed in oder to make meaningful advances.
%------------------------------------------------------------------------
\section{Notation and terminology}\label{NotTerm}\noindent
%------------------------------------------------------------------------
A \emph{minimal system} in a Banach (or quasi-Banach) space $\XX$ over  the real or complex field $\FF$ will be a sequence
 $\XB=(\xx_n)_{n=1}^\infty$  in $\XX$ for which there is a sequence $\XB^*=(\xx_n^*)_{n=1}^\infty$ in $\XX^*$  such that  $\xx_n^*(\xx_n)=\delta_{n,k}$ for all $k$, $n\in\NN$. If $\XB$ is \emph{complete}, i.e., its closed linear span 
 \[
 [\XB]=[\xx_n \colon n\in\NN]
 \]
 is the entire space $\XX$, then $\XB^*$ is unique, and we call it the \emph{dual minimal system} of $\XX$. In this case, we can associate the biorthogonal system $(\xx_n,\xx_n^*)_{n=1}^\infty$  to $\XB$. Also, we can define for each $A\subseteq\NN$ finite the \emph{coordinate projection} on $A$ relative to $\XB$ as 
 \[
 S_A\colon \XX \to \XX, \quad f\mapsto \sum_{n\in A} \xx_n^*(f) \, \xx_n.
 \] 
 A finite set $A\subseteq\NN$ is said to be a \emph{greedy set} of $f\in\XX$ with respect to the complete minimal system $\XB$ if 
\[
\abs{\xx_n^*(f)}\ge \abs{\xx_k^*(f)}, \quad n\in A, \quad k \in\NN\setminus A,
\]
in which case $S_A(f)$ is said to be a \emph{greedy projection}. The map
 \[
 \Ft\colon \XX \to \FF^\NN, \quad f\mapsto \xx_n^*(f)
 \]
 will be called the \emph{coefficient transform} with respect to $\XB$. The support of a function (or signal) $f\in\XX$  with respect to $\XB$ is the set
\[
\supp(f)=\{n\in\NN \colon \xx_n^*(f)\not=0\}.
\]

A sequence  $\XB=(\xx_n)_{n=1}^\infty$  in $\XX$ is said to be a \emph{Schauder basis} if for each $f\in\XX$ there is a unique sequence 
\[
\alpha(f)=(a_n)_{n=1}^\infty
\]
in $\FF$ such that $f=\sum_{n=1}^\infty a_n \, \xx_n$. If this series converges unconditionality, $\XB$ is said to be an \emph{unconditonal basis}. It is known that $\XB$ is an unconditional basis if and only if it is a complete minimal system for which the coordinate projections are uniformly bounded. If $C\in [1,\infty)$ is such that 
\[
\norm{S_A} \le C, \quad A\subseteq \NN, \, \abs{A}<\infty,
\]
we say that $\XB$ is $C$-unconditional. Furthermore, unconditional bases satisfy the estimate
\begin{equation}\label{eq:LatUnc}
\norm{\sum_{n=1}^\infty a_n\, \xx_n} \le C \norm{\sum_{n=1}^\infty b_n\, \xx_n}, \quad \abs{a_n}\le \abs{b_n}.
\end{equation}
If \eqref{eq:LatUnc} holds for a given constant $C$, we say that the basis is $C$-\emph{lattice unconditional}. 
A sequence $\XB$ is a Schauder basis for $\XX$ if and only if it is a complete minimal system for which \emph{partial sum projections} 
 \[
 S_{\enbrace{n\in \NN \colon n\le m}}, \quad m\in\NN,
 \]
 are uniformly bounded; besides, $\Ft(f)=\alpha(f)$ for all $f\in\XX$.
 
Two sequences $(\xx_n)_{n=1}^\infty$ and $(\yy_n)_{n=1}^\infty$ in $\XX$ are said to be \emph{equivalent} if there is an isomorphism $T\colon [\XB] \to [\YB]$ such that $T(\xx_n)=\yy_n$ for all $n\in\NN$. If $C\in [1,\infty)$ is such that 
\[
\max\enbrace{\norm{T}, \norm{T^{-1}}} \le C
\]
we say that $\XB$ and $\YB$ are $C$-equivalent. A \emph{symmetric basis} will be a Schauder basis equivalent to all its permutations, and \emph{subsymmetric basis} will be an unconditional basis $\XB=(\xx_n)_{n=1}^\infty$  equivalent to $(\xx_{\varphi(n)})_{n=1}^\infty$ for all increasing maps  $\varphi\colon\NN\to\NN$. If $\XB=(\xx_n)_{n=1}^\infty$ is a symmetric (resp., subsymmetric) basis then there is a constant $C$ such that 
\begin{enumerate}[label=(S)]
\item\label{it:S} $\XB$ is $C$-equivalent to $(\xx_{\varphi(n)})_{n=1}^\infty$ 
\end{enumerate}
for every one-to-one (resp., increasing) map $\varphi\colon\NN\to\NN$. Any symmetric basis is unconditional, hence subsymmetric. If a sequence $\XB$ is $C$-lattice unconditional and \ref{it:S} holds for every 
one-to-one (resp., increasing) map $\varphi\colon\NN \to\NN$, we say that it is a $C$-symmetric  (resp., $C$-subsymmetric) basis. 

If the constant $C$ involved  in the characterization of unconditonality, lattice unconditionality, symmetry, or subsymmetry is $1$, we say that the corresponding property holds isometrically. Since these properties  are linear, that is, their definitions only involve the boundedness of certain linear operators, any unconditional (resp., symmetric or subsymmetric) basis becomes isometrically lattice unconditional  (resp., symmetric or subsymmetric) under a suitable renorming of the space (see \cite{Ansorena2018} for the slightly more subtle handling of subsymmetric bases). 

A \emph{symmetric sequence space} will be quasi-Banach space $\Sym\subseteq\FF^\NN$ for which the unit vector system is an isometrically symmetric basis of its closed linear span in $\Sym$.
In greedy approximation, several nonlinear forms of symmetry naturally appear. We next introduce some terminology in order to properly define them.  Put
\[
\EE=\{\lambda\in\FF \colon \abs{\lambda}=1\}, \quad\text{and}\quad  \DD=\{\lambda\in\FF \colon \abs{\lambda}\le1\}.
\]
Let $\XB=(\xx_n)_{n=1}^\infty$ be a sequence in a  space $\XX$. Given  $A\subseteq\NN$ finite and $\varepsilon=(\varepsilon_n)_{n\in A}\in\EE^A$, we set
\[
\Ind_{\varepsilon,A}=\sum_{n\in A} \varepsilon_n \, \xx_n, \quad \Ind_{A}=\sum_{n\in A}  \xx_n.
\]
The sequence $\XB$ is said to be $C$-\emph{symmetric for largest coefficients} (SLC for short), $1\le C<\infty$, if 
\begin{equation}\label{eq:SLC}
\norm{ \Ind_{\varepsilon,A} + \sum_{n\in E} a_n\, \xx_n} \le C \norm{ \Ind_{\delta,B} + \sum_{n\in E} a_n\, \xx_n}
\end{equation}
for all pairwise disjoint finite subsets $A$, $B$ and  $E$ of $\NN$, all $\varepsilon\in\EE^A$, all $\delta\in\EE^B$, and all $(a_n)_{n\in E}$ in $\DD$.

If \eqref{eq:SLC} holds in the case when $E=\emptyset$, we say that the sequence $\XB$ is $C$-superdemocratic. If  \eqref{eq:SLC} holds under the aditional restriction that $\varepsilon$ and $\delta$ are constant, we say that  $\XB$ is $C$-democratic. In all cases, if the constant $C$ is irrelevant we drop it from the notation, and if the condition holds with $C=1$, we say that it holds isometrically. Isometric symmetry for largest coefficients was originally named Property~(A) \cite{AW2006}.

A  sequence is superdemocratic if anf only if it is democratic and there is a constant $C$ such that 
\begin{equation}\label{eq:UCC}
\norm{\Ind_{\varepsilon,A}} \le C \norm{\Ind_{\varepsilon,B}}, \quad A\subseteq B\subseteq\NN, \abs{B}<\infty, \, \varepsilon\in\EE^B.
\end{equation}
Sequences that satisfy \eqref{eq:UCC} are called \emph{unconditional for constant coefficients} (UCC for short).

Superdemocracy can be characterized in terms of the  \emph{upper democracy function}, also known as \emph{fundamental function}, $\usdf$ and the \emph{lower democracy function} $\lsdf$ of the sequence. Namely, if for $m\in\NN$ we set
\[
\usdf(m) =\sup_{\substack{\abs{A}\le m \\ \varepsilon\in\EE^A}} \norm{\Ind_{\varepsilon, A}},\quad
\lsdf(m) =\inf_{\substack{\abs{A}\ge m \\ \varepsilon\in\EE^A}} \norm{\Ind_{\varepsilon, A}},
\]
then $\XB$ is $C$-superdemocratic if and only if $\usdf(m)\le C \lsdf(m)$ for $m\in\NN$.

Let $\usdf^*$ denote the fundamental function of the dual system $\XB^*$ of  complete minimal system $\XB$. If
\begin{equation}\label{eq:bidem}
\usdf(m) \lsdf^*(m) \le C m, \quad m\in\NN.
\end{equation}
for some constant $C$, then both $\XB$ and $\XB^*$ are $C$-super-democratic. We call  $C$-\emph{bidemocratic} those complete minimal systems that satisfy \eqref{eq:bidem} for a given constant $C$.

Finally, we say that a Banach space is squeezed between the symmetric sequence spaces $\Sym_1$ and $\Sym_2$ via a complete minimal system $\XB=(\xx_n)_{n=1}^\infty$ if the \emph{series transform}
\[
(a_n)_{n=1}^\infty \mapsto \sum_{n=1}^\infty a_n \, \xx_n
\] 
defines a bounded operator from $\Sym_1$ into $\XX$ and the coefficient transform is a bounded operator from $\XX$ into $\Sym_2$. A complete minimal system is said to be \emph{squeeze symmetric} if it is squeezed between two symmetric sequence spaces that are close to each other in the sense that the fundamental functions of their unit vector systems are equivalent. 

Bidemocratic complete minimal systems are squeeze symmetric. In turn, squeeze symmetry  is stronger than symmetry for largest coefficients.

Any democratic sequence $(\xx_n)_{n=1}^\infty$ is \emph{semi-normalized}, that is,
\[
\inf_{n} \norm{\xx_n}>0, \quad \sup_{n} \norm{\xx_n}<\infty.
\]
In the case when $\XB$  is a Schauder basis, its associated  biorthogonal system $(\xx_n,\xx_n^*)_{n=1}^\infty$ is bounded, i.e.,
\begin{equation}\label{eq:Mbounded}
\sup_n  \norm{\xx_n}  \norm{\xx_n^*}<\infty.
\end{equation}
Since any minimal system becomes normalized under rescaling, and semi-normalization is preserved by equivalence, it is natural to assume the minimal systems we deal with to be semi-normalized. In turn, \eqref{eq:Mbounded} is a feature of the complete minimal system $\XB$ that is convenient to impose in order to implement the \emph{thresholding greedy algorithm} (TGA for short) with respect to it. Note that if $\XB$ is a semi-normalized complete minimal system whose associated biorthogonal system $(\xx_n,\xx_n^*)_{n=1}^\infty$ is bounded, then
\[
\sup_n \norm{\xx_n^*}<\infty.
\]
 Hence, the \emph{coefficient transform} maps $\XX$ into $c_0$. Consequently, for each $f\in\XX$ and $m\in\NN$ there is a (not necessarily unique) greedy set $A$ of $f$ with $\abs{A}=m$. Let $A_m(f)$ denote the one set for which $\max A$ is minimal. The TGA is the sequence $(\GG_m)_{m=1}^\infty$ of nonlinear operators given by
 \[
 \GG_m\colon\XX\to \XX, \quad f\mapsto S_{A_m(f)}.
 \]
 
 Given $f\in\XX$, an arbitrary small perturbation of $f$ yields a  signal $g\in\XX$ for which the greedy sets of any cardinality are unique. This observation leads to the paradigm that any functional property of the mapping $f\mapsto A_m(f)$, $f\in\XX$, yields a property of the  mapping
 \[
 f\mapsto \enbrace{S_A(f) \colon A \mbox{ greedy set of f},\, \abs{A}=m}, \quad f\in\XX.
 \]
 
For simplicity, throughout this paper we will use the term \emph{basis} to refer to a semi-normalized complete minimal system whose associated biorthogonal system is bounded. A \emph{basic sequence} will be a sequence that is a basis of its closed linear span.
 
%---------------------------------------------------
\section{Greedy bases from an isometric point of view}\noindent
%---------------------------------------------------
A basis $\XB$ of a Banach space $\XX$ is said to be \emph{greedy} if the TGA relative to the basis provides optimal sparse approximations, that is,  there $1\le G<\infty$ such that
\begin{equation}\label{defgreedy}
\norm{f-S_A(f)}\le G \norm{f- g},
\end{equation}
whenever  $g$ is a linear combination of $m$ vectors from $\XB$ and $A$ is a greedy set of  $f\in\XX$ with $\abs{A}=m$. If  inequality \eqref{defgreedy} holds for a certain constant $G\ge 1$, we say that $\XB$ is $G$-greedy.

 Konyagin and Temlyakov \cite{KoTe1999} proved that $\XB$ is a greedy basis if and only if it is  demotratic and unconditional. Quantitatively, making the most of their methods, it can be proved that $G$-greedy bases are $G$-unconditional and $G$-democratic; and 
 conversely, if $\XB$ is $C$-lattice unconditional and $D$-democratic, then it is $G$-greedy, where $G=C(1+D)$. The latter estimate gives that isometrically democratic and isometrically lattice-unconditional bases are just $2$-greedy, so it is useless for studying isometrically greedy bases. 
 
 The motivation behind this Section lies in the analysis of the optimality of the thresholding greedy algorithm relative to bases in Banach spaces. This optimality is reflected in the sharpness of the constants that appear in the definitions of the different types of greedy-like bases. What justifies studying the ``isometric'' case in general is the fact that various approximation algorithms converge trivially when some appropriate constant is 1. 
 
 The first movers in this direction were Albiac and Wojtaszczyk, who in \cite{AW2006} characterized isometrically greedy bases using symmetry for largest coefficients instead of democracy. In fact, any $G$-greedy basis is $G$-unconditional and $G$-SSL; and any $C$-unconditional $D$-symmetric for largest coefficients basis is $CD$-greedy (see \cite{AlbiacAnsorena2017b}*{Remark 3.8}). In particular, a basis is isometrically greedy if and only if it is isometrically unconditional and has Property~(A). These estimates opened the door to study the following general question.
 
\begin{question}\label{question:1GUT}
Given a greedy basis  of a Banach space $\XX$, does it become $1$-greedy under a suitable renorming of the space?
\end{question}

Any Banach space with a greedy basis can be renormed so that the basis becomes isometrically unconditional while the SSL constant does not increase. So, Question~\ref{question:1GUT} reduces to the problem of  finding a renorming of $\XX$ so that the basis becomes isometrically symmetric for largest coefficients.

Isometrically symmetric bases are isometrically greedy, thus Banach spaces with a symmetric basis can be renormed so that the basis becomes  isometrically greedy. Despite the fact that there are isometrically subsymmetric bases that are not isometrically greedy \cite{AW2006}, the answer to Question~\ref{question:1GUT} seems to be positive for all the subsymmetric bases found in the literature. For instance, this is the case with Garling sequence spaces,  modelled after an example of Garling from \cite{Garling1968} (cf.  \cite{AAW2018b}).

\begin{problem}\label{question:1SS}
Does any subsymmetric basis of any Banach space become $1$-greedy under a suitable renorming of the space?
\end{problem}

There are isometrically greedy bases that are not subsymmetric \cite{DOSZ2011}*{Theorem 6.9}. So, Question~\ref{question:1GUT} also makes sense for greedy bases that are not subsymmetric. The Haar system in $L_p:=L_p([0,1])$, $p\in(1,2)\cup(2,\infty)$, is probably the most important of example of such bases (see \cite{Temlyakov1998c}).

\begin{problem}\label{question:1Haar}
Let $1<p<\infty$, $p\not=2$. Is there a renorming of $L_p$ so that  $L_p$-normalized Haar system becomes $1$-greedy?
\end{problem}

Notice that the Haar system in $L_p$ for $1<p<\infty$ is bidemocratic, and so are the subsymmetric bases of any Banach space. Thus answering Question~\ref{question:1GUT} in the positive for \emph{bidemocratic} greedy bases would also answer in the positive Problems~\ref{question:1SS} and ~\ref{question:1Haar}. In this regard we mention the following approximation to a solution  of Problem~\ref{question:1Haar}.
\begin{theorem}[\cite{DKOSZ2014}*{Proposition 1.1}]\label{thm:1egreedy}
Let $\XB$ be a  bidemocratic greedy basis of a Banach space $\XX$ and let $C>1$. Then there is a renorming of $\XX$ so that $\XB$ becomes   $C$-greedy.
\end{theorem}

A  more specific question than Question~\ref{question:1GUT}    that still covers  Problem~\ref{question:1Haar} is whether  it admits a positive answer for spaces with nontrivial type. In fact, the fundamental function $\usdf=(s_m)_{m=1}^\infty$ of any superdemocratic basis of a Banach space with nontrivial type has the URP, that is,
\[
s_{rm}\le \frac{r}{2}  s_m \quad m\in\NN,
\]
for some $r\in\NN$ (see \cite{DKKT2003}*{Proof of Proposition 4.1}). Besides, greedy bases, or, more generally, squeeze symmetric bases whose fundamental function has the URP, are  bidemocratic (see \cite{AABW2021}*{Lemma 9.8 and Proposition 10.17(iii)}). In this regard, it is natural to wonder whether the existence of a lattice structure on $\XX$ could aid to obtain a positive answer to Question~\ref{question:1GUT}.

\begin{problem}
Let $\XB$ be a greedy basis of a superreflexive Banach lattice $\XX$.  Is there a renorming of $\XX$ so that $\XB$ becomes $1$-greedy?
\end{problem}

The fundamental function of  the Haar system in $L_p$, $1<p<\infty$, grows as 
 $(m^{1/p})_{m=1}^\infty$ \cite{Temlyakov1998}, so we could also address Problem~\ref{question:1Haar} by focussing on the study of greedy basis whose fundamental functions grow as  $(m^{\alpha})_{m=1}^\infty$ for some $0<\alpha<1$. We point out that the  answer to Question~\ref{question:1GUT} is negative for  for greedy bases other than the  canonical basis of $\ell_1$,  whose fundamental function grows as $(m)_{m=1}^\infty$ (see \cite{DOSZ2011}*{Corollary 5}).

%---------------------------------------------------
\section{Isometric almost greediness}\noindent
%---------------------------------------------------
A basis $\XB$ of a Banach space is said to be \emph{almost greedy} if the TGA provides optimal approximations by means of coordinate projections, that is, there is a constant $1 \le G$ ($G$-almost greedy) such that
\begin{equation*}
\norm{f-S_A(f)}\le G \norm{f- S_B(f)}
\end{equation*}
whenever $A$ is a greedy set of $f$ and $B\subseteq \NN$ satisfies $\abs{A}=\abs{B}$. 
In this Section we are concerned about finding the optimal almost greediness constant.

\begin{question}\label{qt:1AG}
Given an almost greedy basis $\XB$ of a Banach space $\XX$, is there a renorming of $\XX$ so that $\XB$ becomes isometrically almost greedy or, at least, $C$-almost greedy with $C$ arbitrarily close to $1$?
\end{question}

Dilworth et al.\ provided a  characterization of almost greedy bases that runs parallel to that of  greedy ones. For that they used a weaker form of uncondtionality (namely, quasi-greediness) which was introduced in \cite{KoTe1999}. We recall that basis $\XB$ is said to be $C$-\emph{quasi-greedy} for some $1\le C$ if $\norm{g} \le C \norm{f}$ for all $f\in\XX$ and all greedy projections $g$ of $f$. Of course, $\XB$ is quasi-greedy if and only if there is a (possibly different) constant $C$ such that 
\[
\norm{f-g}\le C \norm{f}
\]
for all $f\in\XX$ and all greedy projections $g$ of $f$. If the latter inequality holds for some constant $C\ge 1$, we say that $\XB$ is $C$-\emph{suppression quasi-greedy}.

\begin{theorem}[\cite{DKKT2003}*{Theorem 3.3}]\label{thm:AGQGD}
Let $\XB$ be a basis of a Banach space. Then $\XB$ is almost greedy if and only if it is democratic and quasi-greedy.
\end{theorem}

However,  Theorem~\ref{thm:AGQGD}
does not give quantitative estimates that could aid to address Question~\ref{qt:1AG}. Notwithstanding, it is known \cite{AlbiacAnsorena2017b} that any $G$-almost greedy basis is  $G$-symmetric for largest coefficients and $G$-suppression quasi-greedy; and conversely, $D$-symmetric for largest coefficient and $C$-suppression quasi-greedy bases are $CD$-almost greedy.  Hence, the almost greedy constant is close to one if and only if both the SLC  and the suppression quasi-greedy constants are. 

Being almost greedy is a weaker condition than being greedy.  Taking this into consideration we could expect that  obtaining renormings that improve the almost greedy constant would be easier than improving the greedy constant by renorming  the space. However, in order to tackle  Question~\ref{qt:1AG} we will have to face a new obstruction: since quasi-greediness is not a linear condition, we should develop techniques for improving the suppression quasi-greedy constant. Another general question arises.

\begin{question}\label{qt:1QG}
Given a quasi-greedy basis $\XB$ of a Banach space $\XX$, is there a renorming of $\XX$ so that $\XB$ becomes $1$-suppression quasi-greedy or, at least, $C$-suppression quasi-greedy with $C$ arbitrarily close to $1$?
\end{question}

Since $1$-quasi greedy bases are $1$-suppression unconditional \cite{AlbiacAnsorena2016c} there is no point in replacing the suppression quasi-greedy constant with the quasi-greedy constant in the isometric part of Question~\ref{qt:1QG}.
In contrast, the following problem seems to be unsolved.

\begin{problem}
Is there a Banach space with a conditional $1$-suppression quasi-greedy basis?
\end{problem}

 Turning back to almost greedy bases, we raise the problem whether it is possible to obtain the almost greedy version of Theorem~\ref{thm:1egreedy}.
 \begin{problem}\label{pbm:AIAG} 
Given a bidemocratic almost greedy basis $\XB$ of a Banach space  $\XX$ and $C>1$, is there a renorming of $\XX$ with respect to which $\XB$ is $C$-almost greedy?
 \end{problem}
 
When addressing Question~\ref{qt:1QG} in the isometric case, we should take into account that a basis is $1$-almost greedy if and only if it has Property~(A) \cite{AlbiacAnsorena2017b}, so the problem of improving the suppression unconditional constant disappears in this case, and we must take only care of improving the SLC constant. Besides, once it is shown that isometric symmetry for largest coefficients implies quasi-greediness, we should clarify whether it also implies a stronger condition.

A recent construction from \cite{AABCO2024} showed  the existence of a renorming of $\ell_1$ with respect to which the unit vector system has Property~(A) but it is not $1$-unconditional. However, this result does not answer the main question concerning Property~(A).

\begin{problem}
Does Property~(A) imply unconditionality?
\end{problem}

It is even unknown whether there is a constant $C$ so that any basis with Property~(A) is $C$-unconditional. Note that Property~(A) implies $1$-suppression quasi-greediness. It also unknown whether $1$-suppression quasi-greediness implies unconditionality. 
%---------------------------------------------------
\section{Squeezing spaces between Lorentz spaces}\noindent
%---------------------------------------------------
In this section we will need the dual property of the URP. We say that a sequence $(s_m)_{m=1}^\infty$ in $(0,\infty)$ has the \emph{lower regularity property} (LRP for short) if there is an integer $r\ge 2$ such
\begin{equation*}
2 s_m \le s_{rm}, \quad m\in\NN.
\end{equation*}

Given $0<q\le \infty$ and a nondecreasing  sequence  $\prim=(s_m)_{m=1}^\infty$, we adopt the convention that $s_0=0$, and we define the Lorentz sequence space
\[
d_q(\prim)=\enbrace{ f\in c_0 \colon \Vert f\Vert_{q,\ww}:=\left( \sum_{n=1}^\infty ( b_n\, s_n)^q \frac{s_n-s_{n-1}}{s_n}\right)^{1/q}<\infty},
\]
 where $(b_n)_{n=1}^\infty$ is the nonicreasing rearrangement of $\abs{f}$ and we use  the usual modification if $q=\infty$. If $\prim$ is doubling, that is
\[
 \sup_m  \frac{s_m}{s_{\ceil{m/2}}} <\infty
\]
then
$d_q(\prim)$ is a quasi-Banach space. In fact, $d_q(\prim)$ is a Banach space provided that $1\le q<\infty$. However, $d_\infty(\prim)$ is a Banach space if and only if $\prim$ has the URP. Besides, this spaces are reflexive if and only $1<q<\infty$ and $\prim$ has the URP, and superreflexive if and only if $1<q<\infty$ and $\prim$ has both the URP and the LRP (see \cite{ABW2023}).

For a fixed doubling sequence $\prim$, $(d_q(\prim))_{q>0}$ is an increasing family of symmetric sequence spaces whose fundamental function grows as $\prim$. If $\prim=(n^{1/p})_{m=1}^\infty$ for some $0<p<\infty$, then $d_q(\prim)$ is the classical sequence Lorentz space $\ell_{p,q}$.

Given a basis $(\xx_n)_{n=1}^\infty$ with dual basis $(\xx_n^*)_{n=1}^\infty$ we define for each $f\in\XX$
\[
\varepsilon(f)=(\sgn(\xx_n^*(f)))_{n=1}^\infty\in\EE^\NN,
\]
where $\sgn(0)=1$ and $\sgn(\lambda)=\lambda/\abs{\lambda}$ otherwise. 

The proof of Theorem~\ref{thm:AGQGD} heavily depends on  proving first that if a basis $\XB$ is quasi-greedy then the \emph{restricted truncation operators}
\[
\Rt_m \colon \XX \to \XX, \quad  f\mapsto \min_{n\in A_m(f)}\abs{\xx_n^*(f)}\Ind_{\varepsilon(f),A_m(f)}, \quad m\in\NN,
\]
are uniformly bounded. This condition, coined as \emph{truncation quasi-greediness} (TQG for short) by the authors of \cite{AABBL2022}, is equivalent to the existence of a constant $C$ so that 
\[
 \min_{n\in A}\abs{\xx_n^*(f)} \norm{\Ind_{\varepsilon(f),A}} \le C \norm{f}
 \]
for all $f\in\XX$ and all greedy sets $A$ of $f$. 

A basis $\XB$ is truncation quasi-greedy if and only if the coefficient transform maps $\XX$ into $d_\infty(\lsdf)$. Moreover, for any basis $\XB$ the series transform defines a bounded operator from $d_1(\usdf)$ into $\XX$  (see \cite{AABW2021}*{Section 9}). These results yield that a basis is squeeze symmetric if and only if it is truncation quasi-greedy  and democratic; consequently, almost greedy bases are squeeze symmetric (cf. \cite{AlbiacAnsorena2016}). We also infer that if $\XB$ is  squeeze symmetric then, $\XX$ is squeezed between $d_1(\usdf)$ and  $d_\infty(\usdf)$ via $\XB$. 

While $d_1(\usdf)$ is a Banach space, $d_\infty(\usdf)$ could be nonlocally convex. In fact, if  $\usdf(m)\approx m$ and the coefficient transform maps $\XX$ into a symmetric Banach space then the basis $\XB$ is equivalent to the unit vector vector system of $\ell_1$. To ensure that an almost greedy basis can be squeezed between two Banach spaces we must assume additional conditions such as $\XX$ having nontrivial type so that $\usdf$ verifies the URP.  

When we want to sandwich a Banach space $\XX$  between two symmetric spaces $\mathbb S_{1}$, $\mathbb S_{2}$ that witness that the basis of $\XX$ is squeeze-symmetric, in general we cannot guarantee that $\mathbb S_{1}$ and $\mathbb S_{2}$ retain all the features of $\XX$. It may happen that $\XX$ is locally convex, for instance, but that  $\mathbb S_{2}$ is not. Or that $\XX$ has nontrivial type but we loose that feature in one of the squeezing spaces.
 When dealing with superreflexive spaces this inconvenience disappears. In fact, if $\XB$ is a quasi-greedy basis of a superreflexive Banach space then there is $r>1$ such that the series transform defines a bounded operator from $d_r(\usdf)$ into $\XX$ (see \cite{ABW2023}). Besides, if the basis is almost greedy, a duality argument yields that the coefficient transform is a bounded operator from $\XX$ into $d_q(\usdf)$, whence $\XX$ is squeezed between the superreflexive spaces $d_r(\usdf)$ and $d_q(\usdf)$. We wonder whether superreflexive spaces with a quasi-greedy basis can be squeezed following a similar pattern.

\begin{problem}\label{prob:QGSR}
Let $\XB$ be a quasi-greedy basis of a  superrflexive Banach space $\XX$. Is there $q<\infty$ such that the coefficient transform is a bounded operator from $\XX$ into $d_q(\lsdf)$?
\end{problem}

Here we point out that the answer to the question in Problem~\ref{prob:QGSR} is positive for semi-normalized unconditonal bases (see \cite{AlbiacAnsorena2022c}*{Theorem 7.3}).

%---------------------------------------------------
\section{The TGA and Elton near unconditionality}\noindent
%---------------------------------------------------
A long standing question in basis theory, which was solved in the negative by Gowers and Maurey \cite{GowersMaurey1993}, asked whether all Banach spaces  contained an unconditional basic sequence.  Bearing in mind Rosenthal's theorem \cite{Rosenthal1974}, which states that any bounded sequence in a Banach space either is equivalent to the canonical $\ell_1$-basis or has a weakly Cauchy subsequence,  the most natural way to look for a positive answer to this question was proving that any semi-normalized weakly null sequence has an unconditional basic sequence. When Maurey and Rosenthal \cite{MaureyRosenthal1977} solved in the negative this question, the problem turned to finding subsequences of weakly null sequence that satisfy weaker forms of unconditionality. 

In this ambience,  Elton \cite{Elton1978} introduced near unconditional bases and proved that any normalized weakly null sequence of a Banach space contains a  nearly unconditionality subsequence. Suppose $\XB$ is a basis of a Banach space. Put
\[
 \Cu=\{f\in \XX  \colon \norm{f}_\infty:=\sup_n \abs{\xx_n^*(f)} \le 1\}.
\]
Given a  number $a\ge 0$ and $f\in \XX$ put
\[
A(a,f):=\{n\in\NN \colon \abs{\xx_n^*(f)}\ge a\}.
\] 

The basis $\XB$ is said to be \emph{nearly unconditional} (NU for short) if for each $a\in(0,1]$ there is a constant $C$ such that

\begin{equation}\label{eq:THF}
\norm{ S_A(f)}\le C \norm{f}, \quad f\in\Cu, \, A\subseteq A(a,f).
\end{equation}
The \emph{unconditionality threshold  function} 
\[
\phi\colon(0,1]\to [1,\infty)
\]
is defined for each $a$ as the smallest value of the constant $C$ in \eqref{eq:THF}. Since a basis is unconditional if and only if $\phi$ is bounded, the  unconditionality threshold  function can be used to measure how far q basis iw from being unconditional. We can also measure this distance by means of the \emph{unconditionality parameters}
\[
\unc_m:=\sup_{\abs{A} \le m} \norm{S_A}, \quad m\in\NN
\]
The latter  way of measuring unconditionality is coarser than the former. In fact,
\[
\unc_m \le 2 \sup_n \norm{\xx_n}  \sup_n \norm{\xx_n^*}  \phi(1/m), \quad m\in\NN.
\]
(see \cite{AAB2023b}*{Lemma 6.1}).

 Analogously to greedy, almost greedy,  and squeeze symmetric bases, there is an unconditionality-like condition which combined with democracy characterizes symmetry for largest coefficients.  This condition is called quasi-greediness for largest coefficients, or QGLC for short. We say that a basis $(\xx_n)_{n=1}^\infty$ is $C$-QGLC if
\[
\norm{\Ind_{\varepsilon,A}} \le C \norm{\Ind_{\varepsilon,A}+f}
\]
for all $A\subseteq \NN$ finite, all $\varepsilon\in\EE^A$, and all $f\in\XX$ with $\norm{\Ft(f)}_\infty \le 1$ and $\supp(f)\cap A=\emptyset$. 

Oddly enough, quasi-greediness for largest coefficients and near unconditionality are the same property seen from  different angles \cite{AAB2023}. Besides, the  unconditionality threshold function of any QGLC basis satisfies
\begin{equation}\label{eq:ThPol}
\phi(a)=C a^{-\delta}, \quad 0<a<1
\end{equation}
for some $C\in[1,\infty)$ and $\delta\in(0,\infty)$. 

Truncation quasi-greedy bases of Banach spaces fulfil a better estimate.  Namely,
\begin{equation}\label{eq:ThLog}
\phi(a)=C ( 1-\log a), \quad 0<a<1,
\end{equation}
for some constant $C$ (see \cite{AAB2023b}*{Theorem 6.5}).

While QGLC bases need not be TQG \cite{AAB2024}, it seems to be unknown whether the estimate  \eqref{eq:ThPol} is optimal for QGLC bases.

\begin{problem}
Is there a nearly unconditional basis whose  threshold unconditionality function does not have a logarithmic growth?
\end{problem}

 Elton's  aforementioned subsequence extraction principle has been improved. In fact, it is known that any semi-normalized weakly null sequence has a truncation quasi-greedy subsequence (see \cites{DOSZ2009,AlbiacAnsorena2024}).
However, solving the quasi-greedy basic sequence problem has proven to be a more elusive task.

\begin{problem}\label{pbm:QBSP}
Does any Banach space has a quasi-greedy basic sequence?
\end{problem}

Any Banach space with an unconditonal spreading model which is nonequivalent to the canonical $c_0$-basis has a quasi-greedy basic sequence (\cite{DOSZ2009}, cf.\cite{AlbiacAnsorena2024}). So, in order to address Problem~\ref{pbm:QBSP} it would suffice to focus on Banach spaces without an unconditional basis and whose unique unconditional spreading model is  the standard $c_0$-basis.

%---------------------------------------------------
\section{Semi-greedy bases}\noindent
%---------------------------------------------------
A basis $\XB$ of a Banach space $\XX$ is said to be \emph{semi-greedy} if there is a constant $C$ such that for all $f\in\XX$ and all greedy sets $A$ of $f$ there is $
 h\in [\xx_n \colon n\in A]$
such that $ \norm{f-h}  \le C \norm{f-g}$  for all $g\in\XX$ with $\abs{\supp(g)} \le \abs{A}$. This condition can be reformulated in terms of the  Chebychev-type greedy algorithm, that assigns to each $f\in\XX$ and $m\in \NN$ a vector $\Ct_m(f)$  that minimizes  $\norm{f-h}$ for  $h \in  [\xx_n \colon n\in A_m(f)]$.

Dilworth et al.\ \cite{DKK2003} proved in 2003  that almost greedy bases are semi-greedy. Twenty years later,  Berasategui and Lassalle \cite{BL2023} proved that the converse also holds.

Semi-greedy bases, as well as all  the types of bases we have considered so far, can be defined analogously in the wider framework of  (not necessarily locally convex) quasi-Banach spaces. Most  results on greedy-like bases that were originally stated and proved in Banach spaces hold for quasi-Banach spaces, although the constants involved could be worse. For instance, this  is the case with the characterizations of greedy, almost greedy, squeeze symmetric, symmetric for largest coefficients, and super-democratic bases mentioned above  (see \cite{AABW2021}).  However, it is unknown whether semi-greedy bases behave  in the same in non-locally convex quasi-Banach spaces.

\begin{question}
Is any semi-greedy basis of any quasi-Banach space almost greedy?
\end{question}

\begin{figure}
\begin{equation*}
\begin{gathered}
\xymatrix{
&\text{Greedy} \ar@<3pt>@{=>}[r]\ar@{=>}[d]&\text{Unconditional} \ar@<3pt>@{==>}[l] \ar@{=>}[d]&\\
\text{Semi-greedy} \ar@<3pt>@{-->}[r]&\text{Almost greedy}\ar@<3pt>@{=>}[l]\ar@<3pt>@{=>}[r]\ar@{=>}[d]&\text{Quasi-greedy}\ar@<3pt>@{==>}[l] \ar@{=>}[d]&\\
\text{Bidemocratic}\ar@{=>}[r]&\text{Squeeze symmetric}\ar@<3pt>@{=>}[r]\ar@{=>}[d]&\text{TQG} \ar@<3pt>@{==>}[l] \ar@{=>}[d]&\\
&\text{SLC}\ar@<3pt>@{=>}[r]\ar@{=>}[d]&\text{QGLC} \ar@<3pt>@{==>}[l] \ar@{=>}[d]&\ar@{<=>}[l]\text{NU}\\
&\text{Super-democratic} \ar@<3pt>@{=>}[r]&\text{UCC} \ar@<3pt>@{==>}[l]&\\
}
\end{gathered}
\end{equation*}
\caption{Scale of greedy-like bases}
 \label{fig:scale}
\end{figure}
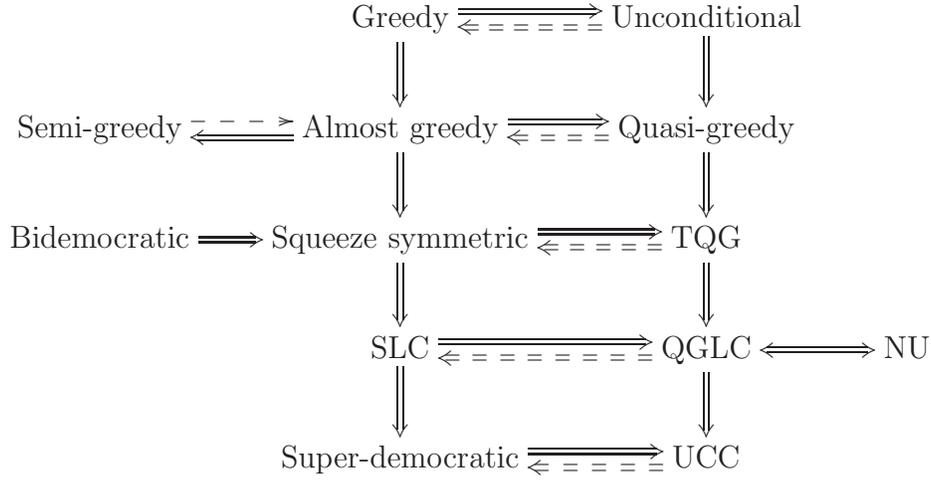

Figure~\ref{fig:scale} represents the relations between the different forms of greedines and unconditionality we have considered. A double arrow means an implication between the two classes of bases involved. A double dashed arrow means that the implication holds under the extra assumption that the basis is democratic. A single  dashed arrow means that the implication holds under the extra assumption that the space is locally convex.

%---------------------------------------------------
\section{Existence of greedy bases}\noindent
%---------------------------------------------------
There are well-known separable Banach spaces, such  as $L_1$, without an unconditional basis, hence without a greedy basis. If a Banach space has a (normalized) unconditional basis we find instances where this basis is, additionally,  democratic thanks to the geometry of the space,  hence that basis ends up being greedy. For example, this is what happens with subsymmetric bases, with the unit vector system of Tsirelson's space, or with the Haar system of both the dyadic Hardy space $H_1$ and   $L_p$ for $1<p<\infty$. On the other hand, we find classical Banach spaces  whose natural basis is unconditional and non-democratic  which might end up having a greedy basis or not. For instance, neither $\ell_p\oplus\ell_q$ nor  the mixed-norm sequence spaces
\[
Z_{p,q}:=\ell_q(\ell_p),
\]
$p$, $q\in[1,\infty]$, $p\not= q$, (if some index is $\infty$ we replace it for $c_0$) have a greedy bases \cites{EdWo1976,Schechtman2014}. However, for $p$ and $q$ in the same range on indices, the Besov spaces 
\[
B_{p,q}=\enpar{\oplus_{n=1}^\infty \ell_p^n}_{\ell_q}
\]
have a greedy basis if and only if  $1<q<\infty$ \cite{DFOS2011}. The original Tsirelson's space does not even have a democratic basis \cite{DKK2003}. 
Note that the space $L_p$ for $p\in(1,2)\cup(2,\infty)$ has a greedy basis but its complemented subspace $\ell_p(\ell_2)$ does not!

Some problems remain open within this research topic. As far as sequence spaces is concerned, the more relevant problem seems to be  settling the greedy basis structure of Nakano spaces.  

Given $\pp=(p_n)_{n=1}^\infty$ in $[1,\infty)$, the Nakano space (a.k.a.\ variable-exponent Lebesgue space) $\ell(\pp)$ consists of all sequences $(a_n)_{n=1}^\infty$ such that 
\[
\sum_{n=1}^\infty \abs{a_n}^{p_n}<\infty.
\]
If $\sup_n p_n=\infty$ the space is not separable, in which case we replace it for its separable part. 

The unit vector system is not a greedy basis for $\ell(\pp)$ unless $\ell(\pp)=\ell_p$ for some $p$ \cite{AADK2016}. If the sequence $\pp$ converges to $1$ or to $\infty$, every complemented unconditional basis sequence of $\ell(\pp)$ is equivalent to a subbasis of the canonical basis \cite{CasKal1998}, and we can infer that $\ell(\pp)$ does not have a greedy  unless $\ell(\pp)=\ell_1$ or $\ell(\pp)=c_0$.

\begin{problem}
Let $\pp=(p_n)_{n=1}^\infty$be a sequence in $[1,\infty)$  with $\lim_n p_n\in(1,\infty)$. Does $\ell(\pp)$ have a greedy basis?
\end{problem}

If $\lim_n p_n=2$ it is even unknown whether $\ell(\pp)$ has a unique unconditional basis (up to rescaling, equivalence and permutation).

As for function spaces it is known \cite{Woj2006} that if $\XX$ is a rearrangement invariant Banach space over $[0,1]$, then the Haar system is a greedy basis for $\XX$ if and only if $\XX=L_p$ for some $1<p<\infty$. The question of whether  rearrangement invariant Banach spaces other that Lebesgue  spaces have a geedy basis naturallly arises. Here, the proximity of the space to $L_1$ or $L_\infty$ may play a role,  so it may be convenient to assume $\XX$ to be superreflexive.
We propose starting to explore this question in the particular case of Orlicz spaces.

\begin{problem}
Let $F\colon [0,1]\to [0,\infty)$ be a $q$ convex and $r$-concave Orlicz function for some $q>1$ and $r<\infty$. Does  $\ell_F$ have a greedy basis?
\end{problem}

The linear structure of  nonlocally convex quasi-Banach spaces is more rigid than the structure of Banach spaces. It is plausible  that any result stating the nonexistence of greedy basis for the spaces of a  family of Banach spaces can be extended to their nonlocally convex relatives (see \cite{AABW2021}*{Section 11.3}). In contrast, extending to the nonlocally convex setting an existence result may require a more careful analysis. Since certain locally convex Besov spaces have a greedy basis, we wonder for the existence of greedy bases for locally nonconvex Besov spaces.
  
\begin{problem}\label{pbm:Besov}
Let $0<p<1<q<\infty$. Does $B_{p,q}$ have a greedy basis?
\end{problem}

If $p$, $q\in(0,1]$,  the spaces $B_{p,q}$ and $Z_{p,q}$ have a unique unconditional basis \cites{KLW1990,AlbiacLeranoz2011} whereas $B_{1,2}$ and $Z_{1,2}$ do not \cite{BCLT1985}. As of today it is unknown whether  the spaces $B_{p,2}$ or $Z_{p,2}$  have a  a unique unconditional basis. Answering in the positive Problem~\ref{pbm:Besov} would yield a normalized unconditional basis $\BB$ of $B_{p,2}$ nonequivalent to any permutation of the canonical basis. Since the unit vector systems of $\ell_p$ and $\ell_2$ are the unique democratic subbases of the the canonical basis $\Zt$ of $Z_{p,2}$, $\BB$ would not be equivalent to any subbasis of $\Zt$. Therefore, the direct sum of $\BB$ and $\Zt$ would be basis of an space isomorphic to $Z_{p,2}$ nonequivalent to any permutation of $\Zt$.

%---------------------------------------------------
\section{Existence of almost greedy bases}\noindent 
%---------------------------------------------------
Once we know that a certain Banach (or quasi-Banach) space does not have a greedy basis, we can ask ourselves whether it has, at least, an almost greedy basis.  To tackle this problem,  Dilworth et al.\ \cite{DKK2003}  invented a method for building conditional almost greedy bases which works for a broad class of spaces. To be precise, applying this method, which we call for short the DKK method, yields conditional almost greedy bases for any Banach space that contains a complemented copy of $\ell_1$, and for any quasi-Banach space that contains a complemented copy of a superreflexive symmetric sequence space (see \cites{DKK2003,AADK2018,AABW2023}). 
The DKK method produces almost greedy bases for $Z_{p,q}$ and $\ell_p\oplus \ell_q$ when $\max\enbrace{p,q}>1$ or $\min\enbrace{p,q}\ge 1$. It also gives almost greedy bases for $B_{p,q}$ when $q>1$ or $\min\enbrace{p,q}\ge 1$. If
$p$, $q\in(0,1)$, then neither $\ell_p\oplus \ell_q$ nor $Z_{p,q}$ nor $B_{p,q}$ have an almost greedy basis unless $p=q$ \cite{AABe2023}.
The following five problems remain open in this direction of research.

\begin{problem}\label{pbm:AGBesov}
Let  $0<p<1$ and let $\XX$ be $\ell_1\oplus\ell_p$,  $Z_{1,p}$, $Z_{p,1}$, $B_{1,p}$ or $B_{1,p}$. Does $\XX$  have an almost greedy basis?
\end{problem}

We point out that apart from the DKK method, which does not seem to work to solve Problem~\ref{pbm:AGBesov}, we do not find in the literature many other techniques that can be used for building almost greedy bases. 

The first example of a conditional quasi-greedy basis, which turned out to be also democratic, was given by Konyagin and Temlyakov \cite{KoTe1999}. Although their technique allows small variations (see \cite{AABW2021}*{Theorem 11.39}), it is far from being suited to solve Problem~\ref{pbm:AGBesov}. Other constructions just give almost greedy bases for particular spaces. Let us mention a couple of results that yield bases for spaces to which the DKK method cannot be applied. The first one, which serves to construct a conditional almost greedy basis for $\ell_p$, $0<p\le 1$, asserts that the Lindenstrauss basis is  conditional and almost greedy in $\ell_p$ (see \cites{AAW2021b,DilworthMitra2001}). The second one says that the Haar system is an almost greedy basic sequence of $\BV(\RR^d)$, $d\ge 2$ \cites{Woj2003,CDVPX1999}.

%---------------------------------------------------
\section{A retrospective look at the role of Schauder bases for implementing the TGA}\noindent 
%---------------------------------------------------
The TGA was initially studied for Schauder bases of Banach spaces. Wojtaszczyk  took the lead towards studying it within the framework of complete minimal systems in quasi-Banach spaces but his initiative did not have many followers.  Since the natural order  in the sequence of the positive integers does not play a significant role in implementing the greedy algorithm, it must be conceded that developing the theory for Schauder basis is somewhat unnatural and limiting.  Besides, the fact of not taking for granted \emph{a priori}  Schauder's condition helps isolate the ingredients  in the proofs of the important results. It could be  argued that minimal systems that are not Schauder bases do not appear naturally, to the extent that the following important problem remains open.

\begin{problem}\label{pbm:SchauderQG}
Is there a quasi-greedy basis  that can not be arranged so that it becomes a Schauder basis?
\end{problem}

If we replace quasi-greediness with truncation quasi-greediness in Problem~\ref{pbm:SchauderQG}, the answer to the question is positive. In fact, there are bidemocratic bases $(\xx_n)_{n=1}^\infty$ such that $(\xx_{\varphi(n)})_{n=1}^\infty$ is not a Schauder basis for any permutation $\varphi$ of $\NN$ (see \cite{AABBL2023}).

%---------------------------------------------------
\section{Dual bases of quasi-greedy bases}\noindent
%---------------------------------------------------
There are quasi-greedy bases whose dual basis are not. Take, for instance, the Lindenstrauss basis of $\ell_1$.
In contrast, the dual basis of any bidemocratic  quasi-greedy basis is quasi-greedy \cite{DKKT2003}*{Corollary 5.5}. Consequently, the dual basis of  an almost greedy of a Banach space with nontrivial type is almost greedy.  Since all known methods for building conditional quasi-greedy bases give almost greedy bases, and the dual basis of normalized unconditional bases is obviously quasi-greedy, no example of a quasi-greedy basis of a Banach space with nontrivial whose dual basis is not quasi-greedy is known. In particular, we pose the following problem.

\begin{problem}
Let $1<p<\infty$, $p\not=2$. Does there exist a quasi-greedy basis of $\ell_p$
whose dual basis in not quasi-greedy?
\end{problem}

Note that any quasi-greedy basis  of $\ell_2$ is almost greedy \cite{Woj2000}*{Theorem 3}, whence its dual basis is quasi-greedy. 

%---------------------------------------------------
\section{Banach envelopes of quasi-greedy bases}\noindent
%---------------------------------------------------
Given a quasi-Banach space $\XX$ there is a pair $(\widehat{\XX},J_\XX)$ that satisfies the universal property associated with the set  of all pairs $(\YY,T)$ consisting of a Banach space $\YY$ and a linear contraction $T\colon\XX\to\YY$. We call $\widehat{\XX}$ the \emph{Banach envelope} of $\XX$ and $J_\XX\colon \XX \to \widehat{\XX}$ the \emph{envelope map} of $\XX$.

 If $\XB$ is a basis of $\XX$, then $\widehat{\XB}:=J_\XX(\XB)$ is a basis of $\widehat{\XX}$ that inherits from $\XB$ all its linear properties. For instance, if $\XB$ is normalized and unconditional in $\XX$ then $\widehat{\XB}$ is (semi-normalized and) unconditional in $\widehat{\XX}$.  
 
 As far as nonlinear properties is concerned, there are instances where $\XB$ is greedy and  $\widehat{\XB}$ is not democratic (see \cite{AABW2021}*{Section 11.7}). Determining whether weaker, nonlinear forms of unconditionality pass to the envelope has proven to be a more elusive task.

\begin{problem}
Let $\XB$ be a quasi-greedy basis of a quasi-Banach space $\XX$. Is  $\widehat{\XB}$ a quasi-greedy basis of $\widehat{\XX}$.
\end{problem}

%---------------------------------------------------------
\section{Weak Chebyshev Greedy Algorithm}
\label{CG}\noindent
%------------------------------------------------------------------------------------------------
So far we have concentrated and discussed  a special case of sparse approximation with respect to a basis and a very specific algorithm to carry out such approximation,  namely the Thresholding Greedy Algorithm. In our  last section we mostly continue to discuss the case of bases but instead of the TGA we consider another greedy-type algorithm, which was 
introduced and studied for sparse approximation with respect to an arbitrary dictionary (see \cite{T15}).

Indeed, in many applications  it is convenient to replace a basis by a more  
general system which may be  redundant, in particular, repetitions are allowed. This latter setting is much  more complicated than the former  (the bases case), however there is a  solid justification of the importance of redundant systems in both theoretical questions and in practical  applications in numerical analysis (see for instance  \cites{S, H, Do3}). The reader can find further discussion of this topic in the books \cites{ST, VTbook, VTsparse, VTbookMA} and the survey papers \cites{T2, T3}. 

In a general setting we will be working in a Banach space $\XX$ with a redundant system of elements that is called a dictionary $\cD$.  Recall that  a set of elements (functions) $\cD$ from $\XX$ is a dictionary  if each $g\in \cD$ has norm   one ($\|g\|=1$),
and the closure of $\cD$ is $\XX$.

A  signal (or function) $f\in \XX$ is said to be \emph{$m$-sparse with respect to $\cD$} if it admits a representation $f=\sum_{i=1}^m c_i \, \xg_i$ with  $\xg_i\in \cD$, $i=1$, \dots, $m$. The set of all $m$-sparse elements is denoted by $\Sigma_m(\cD)$. 

For a given finction $f_0\in \XX$, the error of best $m$-term approximation is given by
\[
\sigma_m(f_0,\cD) := \inf_{g\in\Sigma_m(\cD)} \|f_0-g\|.
\]

In a broad sense, we are interested in the following fundamental question of sparse approximation with redundancy.

\begin{question}\label{qt:Tem1}
 How to design a practical algorithm relative to a dictionary that builds sparse approximations comparable (in the sense of error) to best $m$-term approximations? 
\end{question}

Of course, this is too big of a question which needs to be  tackled in specific situations, but the pattern is the same in all of them, namely we introduce and study an approximation method given by a sequence of maps (an algorithm) $\mathcal A= (\mathcal A_{m})_{m=1}^{\infty}$ relative to a dictionary $\mathcal D$ in a Banach space $\XX$. That is, $\mathcal A_{m}(f)$ belongs to $\Sigma_m(\cD)$ for all $f\in \XX$. Obviously, for any $f\in \XX$ and any $m\in \NN$ we have 
\[
\|f-A_m(f,\cD)\| \ge \sigma_m(f,\cD).
\] 
We are interested in such pairs $(\cD,\mathcal A)$ for which the algorithm $\cA$ provides approximation close to best $m$-term approximation. In order to measure the efficiency of this algorithm we introduce the corresponding definitions. 

\begin{definition}\label{D5.2} We say that $\cD$ is an \emph{almost greedy dictionary} with respect to $\cA$ if there exist  constants $C_1$ and $C_2$ such that for any $f\in \XX$ and $m\in\NN$
\begin{equation*}\|f-\cA_{\lceil C_1m\rceil}(f,\cD)\| \le C_2\sigma_m(f,\cD).
\end{equation*}
\end{definition}

If $\cD$ is an almost greedy dictionary with respect to $\cA$ then $\cA$ gives  almost ideal sparse approximation. It provides  a $\lceil C_1m\rceil$-term approximant as good (up to a constant $C_2$) as the ideal $m$-term approximant for every $f\in \XX$. In the case when $C_1=1$ we call $\cD$ a \emph{greedy dictionary}.

We also need a more general definition. 
\begin{definition} Let $\phi(u)$ be a  function such that 
$\phi(u)\ge 1$.  We say that $\cD$ is a \emph{$\phi$-greedy dictionary with respect to the algorithm $\cA$} if there exists a constant $C_3$ such that for any $f\in \XX$ and $m\in\NN$
\begin{equation*}
\|f-\cA_{\lceil\phi(m)m\rceil}(f,\cD)\| \le C_3\sigma_m(f,\cD).
\end{equation*}
\end{definition}

It was shown in \cite{Tsp} that the Weak Chebyshev Greedy Algorithm, which we define momentarily, is a solution to Question~\ref{qt:Tem1} for a special class of dictionaries.  

For a nonzero element $g\in \XX$, by the Hahn-Banach theorem we can pick up a norming (or peak) functional  for $g$, i.e., an element $F_g\in \XX^{\ast}$ with $ \|F_g\|_{\XX^*} =1$ and such that $F_g(g) =\|g\|_{\XX}.$
Let $f_0\in \XX $ be given. Then for each $m\ge 1$ and any \emph{weakness parameter}  $t\in(0,1]$ we have the following inductive definition.

\begin{itemize}[leftmargin=*]
\item Choose  any element $\varphi_m :=\varphi^{c,t}_m \in \cD$  satisfying
\[
|F_{f_{m-1}}(\varphi_m)| \ge t\sup_{g\in\cD}  | F_{f_{m-1}}(g)|.
\]

\item Put
 $\Phi_m := \Phi^t_m := \spn(\varphi_j \colon 1 \le j \le m)$
and then define $G_m := G_m^{c,t}$ to be the best approximant to $f_0$ from $\Phi_m$.

\item Finally,  let
$
f_m := f^{c,t}_m := f_0-G_m.
$
\end{itemize}

The Weak Chebyshev Greedy Algorithm (WCGA for short) (see \cite{T15}) is a generalization for Banach spaces of the Weak Orthogonal Matching Pursuit (WOMP). In a Hilbert space the WCGA coincides with the WOMP. The WOPM is very popular in signal processing, in particular in compressed sensing. In the case when $t=1$,   the WOMP is called  Orthogonal Matching Pursuit (OMP).

We note that the properties of a given basis in a Banach space with respect to the TGA and with respect to the  WCGA could be very different, both from a qualitative and a quantitative approach. 
To illustrate the differences in the performance of the TGA and WCGA we use  the $d$-variate trigonometric system  $\cT^d$ as an example.

\begin{proposition}[\cite{VTsparse}*{Theorem 2.2.1}] 
 Let $1\le p \le \infty$, and set $h(p) = |1/2-1/p|$. There is a constant $C$ such that
\[
\|\mathcal G_m(f,{\mathcal T}^d)\|_p \le Cm^{h(p)}\|f\|_p, \quad m\in\NN, \, f \in L_p(\mathbb T^d). 
\]
Moreover, the extra factor $m^{h(p)}$ cannot be improved. 
\end{proposition}
 The following inequalities were obtained very recently (see \cite{VT197}). 
 
\begin{proposition}
Let $2\le p< \infty$. For any $f\in L_p$ and for each $m\in \mathbb N$,
\begin{equation}\label{TG2}
\|\mathcal G_m(f,\cT^d)\|_p \le C(p)\|f\|_p^{2/p}\|f\|_{A_1(\cT^d)}^{1-2/p}.
\end{equation}
Let $1\le p\le 2$. For any $f\in L_p$ and for each $m\in \mathbb N$,
\begin{equation}\label{TG3}
\|\mathcal G_m(f,\cT^d)\|_p \le C(p)\|f\|_p^{p/2}\|f\|_{A_1(\cT^d)}^{1-p/2},
\end{equation}
where
$$
\|f\|_{A_1(\cT^d)}:= \sum_{\bm k\in \mathbb Z^d} |\hat f(\bm k)|,\quad \hat f(\bm k) := (2\pi)^{-d}\int_{\mathbb T^d} f(\bm x)e^{-i(\bm k,\bm x)}d\bm x.
$$
\end{proposition}

Clearly, the $\|\cdot\|_{A_1(\cT^d)}$ norm is stronger than the $\|\cdot\|_p$ norm. However, it is important that 
in inequalities (\ref{TG2}) and (\ref{TG3}) the extra factor $C(p)$ does not depend on $m$.

To compare those estimates with the ones we get for the  WCGA we need to consider the real trigonometric system because the  Weak Chebyshev Greedy Algorithm has been well studied mainly for  real Banach spaces. The reader can find some results on the WCGA for complex Banach spaces in \cite{DGHKT}. We will denote by ${\mathcal R}{\mathcal T}$ the real trigonometric system 
\[
\{1,\,\sin(2\pi x) ,\, \cos(2\pi x), \dots, \sin(2\pi n x), \, \cos(2\pi n x),\dots\}
\]
on $[0,1]$ and will let ${\mathcal R}{\mathcal T}_p$  be its normalized version  in $L_p([0,1])$. Let 
\[
{\mathcal R}{\mathcal T}_p^d = {\mathcal R}{\mathcal T}_p\times\cdots\times {\mathcal R}{\mathcal T}_p
\]
be the $d$-variate real trigonometric system.  The following Lebesgue-type inequality for the WCGA was proved in \cite{Tsp}.

\begin{proposition}\label{T1.2} 
 Let $2\le p<\infty$,  $d\in\NN$ and $t\in(0,1]$. There are constants $C_1$ and $C_2$ such that the WCGA with weakness parameter $t$ relative to the $d$-variate trigonometric system $\cD = \cR\cT_p^d$ in $L_p$ gives
\begin{equation}\label{I1.4}
\norm{f_{\ceil{C_1 m\ln (m+1)}}}_p \le C_2\sigma_m(f_0,\cD)_p, \quad  m\in\NN,\,  f_0\in L_p([0,1]^d).
\end{equation}
\end{proposition}

\begin{problem}[See \cite{T2}*{Problem 7.1}] \label{Tem0}Does \eqref{I1.4} hold without the $\ln(m+1)$ factor?

\end{problem}

 Proposition~\ref{T1.2} is the first result on the Lebesgue-type inequalities for the WCGA with respect to the trigonometric system. It  is a first step towards the solution of Problem~\ref{Tem0}, but the problem is still open. 

 The dissimilarities between the TGA and the WCGA can be realized as well from a qualitative point of view. Take for instance, the class of quasi-greedy bases relative  to the TGA, which  is a rather narrow subset of bases close in a certain sense to the set of unconditional bases. The situation is dramatically different for the WCGA. For example, if $\XX$ is uniformly smooth then WCGA converges for each $f\in \XX$ with respect to any dictionary in $\XX$ (see \cite{VTbook}*{Ch.6}).  

There exists a general theory of the Lebesque-type inequalities for the WCGA with respect to dictionaries satisfying certain conditions. The reader can find some of the corresponding results in \cite{Tsp} and \cite{VTbookMA}*{Ch.8}. We only present here some corollaries of those general results in the case of 
bases satisfying certain conditions. For $p\in (1,\infty)$ we use the notation $p'= p/(p-1)$. 

\begin{proposition}\label{CGP1} Let $\XB$ be a uniformly bounded orthogonal system normalized in  $L_p(\Omega)$ for $1< p<\infty$, where $\Omega$ is a bounded domain. Then,  given $0<t\le 1$ there are constants $C_1$ and $C_2$ depending on $p$, $t$ and $\Omega$, such that
 \begin{equation*}
\norm{f_{\ceil{C_1 m\ln (m+1)}}}_p \le C_2\sigma_m(f_0,\XB)_p,\quad  m\in\NN, f_0\in L_p(\Omega), \,  2\le p<\infty,
\end{equation*}
and
 \begin{equation*}
\norm{f_{\ceil{C_1 m^{p'-1}\ln (m+1)}}}_p \le C_2\sigma_m(f_0,\XB)_p,\quad  m\in\NN,\,  f_0\in L_p(\Omega), \,  1< p\le 2.
\end{equation*}
\end{proposition}

\begin{proposition}\label{CGP2} 
 Let $2\le p<\infty$, $d\in\NN$ and $t\in(0,1]$. Then, the normalized $d$-variate Haar basis in $L_p$ satisfies, for some constants $C_1$ and $C_2$,
\begin{equation*}
\norm{f_{\ceil{C_1 m^{2/p'}}}}_p \le C_2\sigma_m(f_0,{\mathcal H}^d_p)_p, \quad m\in\NN, \,  f_0\in L_p([0,1]^d).
\end{equation*}
\end{proposition}

 \begin{proposition}\label{CGP3}
 Let $1<p\le 2$ and $t\in(0,1]$. Then the univariate Haar basis in $L_p$ satisfies
\begin{equation*}
\norm{f_{C_1m}}_p \le C_2\sigma_m(f_0,{\mathcal H}_p)_p, \quad m\in\NN, \,  f_0\in L_p([0,1]),
\end{equation*}
 for some constants $C_1$ and $C_2$.
\end{proposition}

\begin{proposition}\label{CGP4}
  Let $\gamma\in(0,\infty)$, $1<q\le 2$ and $t\in(0,1]$. Let $\XX$ be a superreflexive Banach space whose modulus of smoothness $\rho$ satisfies $\rho(u)\le \gamma u^q$ for all $u>0$. Assume that $\XB$ is a normalized Schauder basis for $\XX$. 
 Then there are constants $C_1$ and $C_2$ such that
  \begin{equation*}
\norm{f_{\ceil{C_1 m^{q'}\ln(m+1)}}} \le C_2\sigma_m(f_0,\XB), \quad m\in\NN, \, f_0\in\XX.
\end{equation*}
\end{proposition}

\begin{proposition}\label{CGP5} 
Let $0<t\le 1$ and $\XB$ be a normalized quasi-greedy basis for $L_p$, $1< p<\infty$.   Then \begin{equation*}
\norm{f_{\ceil{C_1 m^{2(1-1/p)} \ln (m+1)}}}_p \le C_2 \sigma_m(f_0,\XB),\quad  m\in\NN, \,  f_0\in L_p, \, 2\le p<\infty,
\end{equation*}
and
\begin{equation*}
\|f_{\ceil{C_1 m^{p'/2} \ln (m+1)}}\|_p \le C_2\sigma_m(f_0,\XB),\quad  m\in\NN,\, f_0\in L_p, \,  1<p\le 2,
\end{equation*}
for some constants $C_1$ and $C_2$.

\end{proposition}

\begin{proposition}\label{CGP6} Let $0<t \le 1$ and $\XB$ be a normalized uniformly bounded orthogonal quasi-greedy basis for $L_p$, $1< p<\infty$ (for existence of such bases see \cites{DS-BT, N}). Then 
\begin{equation*}
\|f_{\ceil{C_1 m  \ln(\ln(m+3))}}\|_p \le C\sigma_m(f_0,\XB)_p,\quad m\in\NN,\,  f_0\in L_p,\,  2\le p<\infty,
\end{equation*}
and
\begin{equation*}
\norm{f_{\ceil{C_1 m^{p'/2}  \ln(\ln(m+3))}}}_p \le C_2\sigma_m(f_0,\XB)_p,\quad m\in\NN, \, f_0\in L_p, \, 1<p\le 2,
\end{equation*}
for some constants $C_1$ and $C_2$.
\end{proposition}

The approximation method provided by WCGA is not as well studied as the TGA from the point of view of the 
Lebesgue-type inequalities and greedy-type bases. It is a very interesting albeit very difficult area of research. The reader can find recent results in this direction in the paper \cite{DGHKT}. We close with some open problems of this aspect of the theory in the special case when we study sparse approximation with respect to bases  instead of with respect to redundant dictionaries (see \cite{VTbookMA}*{p.448}).

\begin{problem} Characterize almost greedy bases with respect to the WCGA  for the Banach space $L_p$, $1<p<\infty$. 
\end{problem}

\begin{problem}  Is the $d$-variate trigonometric system $\mathcal R\mathcal T^d_p$ an almost greedy basis with respect to the WCGA in $L_p(\mathbb T^d)$, $1<p<\infty$?
\end{problem}

\begin{problem}  Is the univariate Haar basis $\mathcal H_p$ an almost greedy basis with respect to the WCGA in $L_p$, $2<p<\infty$?
\end{problem}

\begin{problem}   Is the $d$-variate Haar basis $\mathcal H_p^d$, $d\ge 2$, an almost greedy basis with respect to the WCGA in $L_p$, $1<p<\infty$?
\end{problem}

\begin{problem}  For each $L_p$, $1<p<\infty$, find the best $\phi$ such that any Schauder basis is a $\phi$-greedy basis with respect to the WCGA. 
\end{problem}

\begin{problem}   For each $L_p$, $1<p<\infty$, find the best $\phi$ such that any unconditional basis is a $\phi$-greedy basis with respect to the WCGA. 
\end{problem}

\begin{problem}  Is there a greedy-type algorithm $\cA$ such that the multivariate Haar system $\cH^d_p$ is an almost greedy basis of $L_p$, $1<p<\infty$, with respect to $\cA$?
\end{problem}

\begin{problem}
In all the above Propositions \ref{CGP1}--\ref{CGP2},  \ref{CGP4}--\ref{CGP6} we do not know the right $\phi$ which makes the corresponding bases $\phi$-greedy bases with respect to the WCGA. 
\end{problem}

%-------------------------------------------------- \bib, bibdiv, biblist are defined by the amsrefs package

%------------------------------------------------
\end{document}